\documentclass[11pt, a4paper]{article}
\usepackage{amssymb}
\usepackage{makeidx}
\usepackage{latexsym}
\usepackage{mathrsfs}
\usepackage{eufrak}
\newtheorem{theorem}{Theorem}

\newtheorem{definition}[theorem]{Definition}

\newtheorem{corollary}[theorem]{Corollary}

\begin{document}

\def\P{\mathrm{P}}
\def\Field{\mathrm{Field}}
\def\Rank{\mathrm{Rank}}
\def\0{\varnothing}
\def\R{\mathrm{R}}
\def\ID{\textsf{I}\Delta_0 + \mathsf{exp}}
\def\TC{\mathrm{TC}}
\def\EA{\mathsf{EA}}
\def\ZFinf{\mbox{\textsf{ZF--Inf}}}
\def\Lex{\mathrm{Lex}}
\def\<{\left \langle}
\def\>{\right \rangle}
\def\L{\mathscr{L}}
\def\Ack{\mathrm{Ack}}
\def\M{\mathfrak{M}}
\def\Num{\mathrm{Num}}

\title{On interpretations of bounded arithmetic\\ and bounded set theory\footnote{Department of Philosophy, University of Bristol, 9 Woodland Road, Bristol. BS8 1TB. \textsf{Richard.Pettigrew@bris.ac.uk}}\ \footnote{Thanks to John Mayberry and Ali Enayat on earlier versions of this paper.}
\author{Richard Pettigrew}
}
\date{\today}
\maketitle

\begin{abstract}
In \cite{kaye2007oio}, Kaye and Wong proved the following result, which they considered to belong to the folklore of mathematical logic.
\begin{theorem}
The first-order theories of Peano arithmetic and Zermelo-Fraenkel set theory with the axiom of infinity negated are bi-interpretable.
%: that is, they are mutually interpretable with interpretations that are inverse to each other.
\end{theorem}
In this note, I describe a theory of sets that is bi-interpretable with the theory of bounded arithmetic $\ID$.  Because of the weakness of this theory of sets, I cannot straightforwardly adapt Kaye and Wong's interpretation of the arithmetic in the set theory.  Instead, I am forced to produce a different interpretation.
%Finally, I build on the work of Enayat, Schmerl, and Visser \cite{enayatforthomo} to show that a particular weaker theory of sets is not bi-interpretable with $\ID$.
\end{abstract}

\vspace{5mm}

\noindent \textsf{Primary Subject:} 03C62\medskip


\vspace{5mm}

\section{Introduction}

In \cite{kaye2007oio}, Kaye and Wong proved the following result, which they considered to belong to the folklore of mathematical logic.
\begin{theorem}
The first-order theories of Peano arithmetic and Zermelo-Fraenkel set theory with the axiom of infinity negated are bi-interpretable:  that is, they are mutually interpretable with interpretations that are inverse to each other.
\end{theorem}
More precisely, they showed that $\mathsf{PA}$ and $\ZFinf^*$ are bi-interpretable, where $\ZFinf$ is obtained from $\mathsf{ZF}$ by negating the Axiom of Infinity, and $\ZFinf^*$ is obtained from $\ZFinf$ by adding an Axiom of Transitive Containment, which says that each set is contained in a transitive set.

In this note, I describe a theory of sets that is bi-interpretable with the bounded arithmetic $\ID$.  Because of the weakness of this theory of sets, I cannot straightforwardly adapt Kaye and Wong's interpretation of the arithmetic in the set theory.  Instead, I am forced to produce a different interpretation.

In Section \ref{interps}, I lay down some notation and definitions that will aid the discussion of interpretations throughout the paper.  In Section \ref{EA}, I describe a theory of sets called $\EA^*$ and, in Section \ref{EAID}, I describe the Ackermann interpretation of this theory in $\ID$.  In Section \ref{ordcardinterps}, I consider a natural interpretation of $\ID$ in $\EA^*$, but note that it is not inverse to the Ackermann interpretation, and, in Section \ref{IDEA}, I describe my alternative to Kaye and Wong's inverse.
%Finally, in Section \ref{omegamodels}, I discuss $\omega$-models of a related set theory, called $\EA^*$.

\section{Interpretations}\label{interps}

Kaye and Wong consider only first-order languages with relation symbols.  However, both $\ID$ and the set theory described in Section \ref{EA} are most naturally formulated in languages that include function symbols.  After all, amongst the axioms of both theories are schema that are indexed by the set of \emph{bounded quantifier formulae} of the language---the Axiom Schema of Induction for $\ID$ and the Axiom Schema of Subset Separation for $\EA^*$.  And it is most natural to stipulate which formulae are to count as bounded by appealing to terms of the language built up using function symbols.  Thus, we consider first-order languages with relation symbols \emph{and} function symbols.  However, this is not essential.  We could formulate both theories in languages that contain only relation symbols.  And, if we were to do this, our interpretability results would still go through using Kaye and Wong's definition of bi-interpretability for languages that only contain relation symbols.

Like Kaye and Wong, we demand that each language consider contains a unary relation symbol $\mathrm{Dom}$ and each theory contains the sentence $\forall x \mathrm{Dom}(x)$.

Suppose $\L$ is such a language.  Then an \emph{$\L$-theory} is a consistent set of $\L$-sentences.  Given a theory $T_1$ in language $\L_1$ and theory $T_2$ in language $\L_2$, an \emph{atomic interpretation mapping of $T_1$ into $T_2$} is a mapping $\mathfrak{i}$ such that 
\begin{enumerate}
\item[i.] For each function symbol $f$ of $\L_1$ and free variables $\vec{x}$, $f(\vec{x})^\mathfrak{i}$ is a term of $\L_2$ in the same free variables, and
\item[ii.] For each relation symbol $R$ of $\L_1$ and free variables $\vec{x}$, $R(\vec{x})^\mathfrak{i}$ is a formula of $\L_2$ in the same free variables.
\end{enumerate}
Given an interpretation mapping $\mathfrak{i} : T_1 \rightarrow T_2$, we can extend it to a \emph{full interpretation mapping} (also called $\mathfrak{i}$), which takes any formula in $\L_1$ to a formula in $\L_2$.  We define $(\neg \varphi(\vec{x}))^\mathfrak{i}$ to be $\neg \varphi(\vec{x})^\mathfrak{i}$, $(\varphi(\vec{x}) \rightarrow \psi(\vec{x}))^\mathfrak{i}$ to be $\varphi(\vec{x})^\mathfrak{i} \rightarrow \psi(\vec{x})^\mathfrak{i}$, and $(\forall y \varphi(\vec{x}, y))^\mathfrak{i}$ to be $\forall y(\mathrm{Dom}(y)^\mathfrak{i} \rightarrow \varphi(\vec{x}, y))$.  Given a full interpretation mapping $\mathfrak{i} : T_1 \rightarrow T_2$, we say that \emph{$\mathfrak{i}$ defines an interpretation of $T_1$ in $T_2$} if
\begin{enumerate}
\item[i.] $T_2 \vdash \exists x \mathrm{Dom}(x)$, and
\item[ii.] For each sentence $\sigma \in T_1$, $T_2 \vdash \sigma^\mathfrak{i}$.
\end{enumerate}

Now we define two kinds of mutual interpretability, the second stronger than the first:
\begin{enumerate}
\item[(1)] We say that \emph{$T_1$ and $T_2$ are mutually interpretable} if there are interpretations $\mathfrak{f} : T_1 \rightarrow T_2$ and $\mathfrak{g} : T_2 \rightarrow T_1$.
\item[(2)] We say that \emph{$T_1$ and $T_2$ are bi-interpretable} if there are interpretations $\mathfrak{f} : T_1 \rightarrow T_2$ and $\mathfrak{g} : T_2 \rightarrow T_1$ \emph{and}
\begin{enumerate}
\item[i.] for every formula $\varphi$ in $\L_1$, $T_1 \vdash \forall \vec{x} ((\varphi(\vec{x})^\mathfrak{f})^\mathfrak{g} \leftrightarrow \varphi(\vec{x}))$ and
\item[ii.] for every formula $\psi$ in $\L_2$, $T_2 \vdash \forall \vec{x} ((\psi(\vec{x})^\mathfrak{g})^\mathfrak{f} \leftrightarrow \psi(\vec{x}))$
\end{enumerate}
\end{enumerate}
Kaye and Wong proved that $\mathsf{PA}$ and $\ZFinf^*$ are bi-interpretable. In his doctoral thesis \cite{homolka1983ea}, Vincent Homolka described a theory of sets called EA, first formulated by John Mayberry, and proved that it is mutually interpretable with $\ID$.\footnote{Mayberry has since written a book on his system \cite{mayberry2000fmt}.  For the state of the art on this theory, see \cite{mayberrymsnna}.  Independently of Homolka's work, Gaifman and Dimitracopoulos \cite{gaifman1982fpa} had described a theory of sets a year earlier, which they dubbed \textsf{EF}, and which is also mutually interpretable with $\ID$.  I restrict my attention to Mayberry's theory as studied by Homolka.} Here, I describe an extension of EA, which I call $\EA^*$, and I prove that $\EA^*$ is bi-interpretable with $\ID$.

%I finish by proving that $\EA^-$, a particular weakening of EA, is \emph{not} bi-interpretable with $\ID$ using techniques developed by Enayat, Schmerl, and Visser \cite{enayatforthomo}.

\section{A bounded theory of finite sets}\label{EA}

In this section, I describe Mayberry's theory, $\EA$.  Essentially, $\EA$ is obtained from $\mathsf{ZF}$ set theory in three steps:  replace the Axiom of Infinity by an axiom that states that every set is Dedekind finite; restrict the Separation and Replacement axiom schema to hold only for bounded quantifier formulae; and add an axiom of transitive closure.  In \cite{mayberry2000fmt}, Mayberry also described an extension of $\EA$ that is obtained by adding an axiom that guarantees, for every set, the existence of the first level of the cumulative hierarchy at which that set occurs: he calls this axiom the \textsf{Weak Hierarchy Principle} (henceforth, \textsf{WHP}).  In this note, I will consider the theory
\[
\EA - \mbox{\textsf{Transitive Closure}} - \mathsf{Replacement} + \mathsf{WHP}.
\]
I will call this theory $\EA^*$.  It is this theory that is bi-interpretable with $\ID$.

$\EA^*$ is a first-order theory.  Like $\ID$, its language contains function symbols:  in $\ID$, these are used to state the restrictions on induction; in $\EA^*$, they are used to state the restrictions on subset separation.  It has one constant symbol, $\0$. (As usual, this is considered as a $0$-place function symbol.) It has three unary function symbols, $\P(\_)$ (power set), $\bigcup(\_)$ (sum set),  and $\mathrm{R}(\_)$ (rank function):  the latter is introduced by the \textsf{Weak Hierarchy Principle}.  It has one binary function symbol: $\{ \_, \_\}$ (pair set).  And, for each bounded quantifier formula $\Phi$, it has the unary function symbol $\{x \in \_ : \Phi(x)\}$ (subset separation for bounded quantifier formula), where a bounded quantifier formula is one in which each occurrence of a quantifier has the form $\forall y(y \in t(\vec{x}) \rightarrow \Phi(\vec{x}, y))$ or $\exists y(y \in t(\vec{x}) \wedge \Phi(\vec{x}, y))$ for some term $t$ of $\EA^*$.

The axioms of $\EA^*$ are \textsf{Extensionality}, \textsf{Pair Set}, \textsf{Sum Set}, \textsf{Power Set}, \textsf{Foundation}, \textsf{Axiom Schema of Subset Separation for Bounded Quantifier Formulae}, \textsf{Dedekind Finiteness}, and the \textsf{Weak Hierarchy Principle}.  We state the latter three precisely.\medskip

\noindent \textsf{Axiom Schema of Subset Separation for Bounded Quantifier Formulae}\medskip

$\forall \vec x \forall y \forall z (z \in \{u \in y : \phi(u, \vec x)\} \equiv z \in y \wedge \phi(z, \vec x))$\medskip

\noindent for each bounded quantifier formula $\phi$.\medskip

\noindent \textsf{Axiom of Dedekind Finiteness}\medskip

$\forall x, f(f : x \rightarrow x \wedge f \mbox{ is one-one} \rightarrow f \mbox{ is onto})$\medskip

\noindent To state Mayberry's \textsf{Weak Hierarchy Principle}, we must say what it means to be a level in the cumulative hierarchy:
\begin{definition}Given a set $S$, we say that \emph{$S$ is a level in the cumulative hierarchy} if there is a (finite) linear ordering $[V_0, ..., V_n]$ such that $V_0, ..., V_n \subseteq S$, $V_0 = \0$, $V_n = S$, and for each $V_k$, $V_{k+1} = \P(V_k)$.\end{definition}
(Note: since $V_0, ..., V_n \subseteq S$, the property of being a level of the cumulative hierarchy is represented by a bounded quantifier formula.)\medskip

\noindent With this in hand, we can state the \textsf{Weak Hierarchy Principle}.\medskip

\noindent \textsf{Weak Hierarchy Principle}\medskip

$\forall x(x \in \R(x) \wedge \R(x) \mbox{ is a level in the cumulative hierarchy } \wedge$

\hspace{10mm} $\forall y(x \in y \wedge y \mbox{ is a level in the cumulative hierarchy} \rightarrow \R(x) \subseteq y))$\medskip

\noindent A remark is in order.  Like Kaye and Wong, I wish to interpret set theory in arithmetic using the interpretation described by Ackermann in \cite{ackermann1937dwd}.  If we are to find an inverse to this interpretation, we must ensure, for every sentence $\sigma$ of our chosen set theory, that our chosen arithmetic proves the Ackermann translation of $\sigma$ only if our set theory proves $\sigma$.  Here are two important examples.  $\mathsf{PA}$ and $\ID$ both prove the \textsf{Axiom of Dedekind Finiteness} and both prove the Weak Hierarchy Principle.  However, neither sentence occurs as an axiom in Kaye and Wong's $\ZFinf^*$.  This is not a problem because, in $\ZFinf$, the \textsf{Weak Hierarchy Principle} is equivalent to Kaye and Wong's $\mathsf{TC}$, which says that each set is contained in a transitive set; and, in \textsf{ZF} without \textsf{Infinity}, the \textsf{Axiom of Dedekind Finiteness} is equivalent to the negation of \textsf{Infinity}.  However, neither of these equivalences hold in the relevant fragments of EA:  \textsf{Transitive Containment} follows from \textsf{Weak Hierarchy} and $\neg$\textsf{Infinity} follows from the \textsf{Axiom of Dedekind Finiteness}, but neither converse holds.  Thus, we must include the full strength of the \textsf{Axiom of Dedekind Finiteness} and the \textsf{Weak Hierarchy Principle} in our axioms.

\section{The Ackermann interpretation}\label{EAID}

As mentioned above, I will exploit Ackermann's interpretation of arithmetic in set theory to interpret $\EA^*$ in $\ID$.  I describe this interpretation in this section; in Section \ref{IDEA}, I describe its inverse.  

The Ackermann interpretation of set theory in arithmetic is based on the following interpretation of the membership relation:
\[
(x \in y)^\mathfrak{a} \mbox{ is } (\exists n < y)(\exists m < 2^x)[y = 2^{x+1}n + 2^x + m]
\]
The right-hand side says that the $x^\mathrm{th}$ bit of $y$ is 1.  Further,
\begin{eqnarray*}
\mathrm{Dom}(x)^\mathfrak{a} & \mbox{ is } & \mathrm{Dom}(x) \\
(x = y)^\mathfrak{a} & \mbox{ is } & (x = y) \\
\0^\mathfrak{a}  & \mbox{ is } & 0.
\end{eqnarray*}
To complete our definition of $\mathfrak{a} : \EA^* \rightarrow \ID$, we must define $\P(x)^\mathfrak{a}$, $\bigcup(x)^\mathfrak{a}$, $\R(x)^\mathfrak{a}$, $\{x, y\}^\mathfrak{a}$, and, for each bounded quantifier formula $\Phi$, $\{x \in y : \Phi(x)\}^\mathfrak{a}$.  These are straightforward to define, if somewhat intricate.  With this in hand, it is equally straightforward to establish that
\begin{theorem}
$\mathfrak{a}$ defines an interpretation of $\EA^*$ in $\ID$.
\end{theorem}
\emph{Proof. } The proofs of Extensionality$^\mathfrak{a}$ and Foundation$^\mathfrak{a}$ are adapted from the well-known proofs of these sentences in $\mathsf{PA}$.  In that case, they are proved by induction.  In $\ID$, we identify bounds for the quantifiers in the induction formulae and proceed as before.  

Dedekind Finiteness$^\mathfrak{a}$ is derived as a consequence of the Ackermann interpretation of what Mayberry calls \emph{One Point Extension Induction}, which is easily seen to be provable by bounded induction in $\ID$ (Theorem 8.3.3 of \cite{mayberry2000fmt}).  One Point Extension Induction says that, for any bounded quantifier formula $\Phi$, we have
\[
[\Phi(\0) \wedge \forall x \forall z( \Phi(x) \rightarrow \Phi(x \cup \{z\}))]\rightarrow \forall x\Phi(x)
\]

The proof concludes by establishing that $\P(x)^\mathfrak{a}$, $\bigcup(x)^\mathfrak{a}$, $\R(x)^\mathfrak{a}$, $\{x, y\}^\mathfrak{a}$, and $\{x \in y : \Phi(x)\}^\mathfrak{a}$ have the properties that the translations of the corresponding axioms require of them.  Details can be found in \cite{homolka1983ea}. \hfill $\Box$

\section{The ordinal and cardinal interpretations}\label{ordcardinterps}

Kaye and Wong note that there is an obvious interpretation of $\mathsf{PA}$ in $\ZFinf^*$, which interprets the arithmetic as ordinal arithmetic.  Thus, let $\mathrm{Ord}$ be the class of von Neumann ordinals, as usual, and define the following relations on this class:  $x +_\mathrm{o} y = z$ (ordinal addition) and $x \times_\mathrm{o} y = z$ (ordinal multiplication).  Then let $\mathfrak{o} : \mathsf{PA} \rightarrow \mbox{ZF--Inf}^*$ be the interpretation mapping defined as follows:
\begin{eqnarray*}
\mathrm{Dom}(x)^\mathfrak{o} & \mbox{ is } & x \in \mathrm{Ord} \\
(x = y)^\mathfrak{o} & \mbox{ is } & (x = y) \\
(x < y)^\mathfrak{o} & \mbox{ is } & x \in y \\
(x + y = z)^\mathfrak{o} & \mbox{ is } & (x +_\mathrm{o} y = z) \\
(x \cdot y = z)^\mathfrak{o} & \mbox{ is } & (x \times_\mathrm{o} y = z) \\
\end{eqnarray*} 
\begin{theorem}
$\mathfrak{o}$ defines an interpretation of $\mathsf{PA}$ in $\ZFinf^*$.
\end{theorem}
As Kaye and Wong point out, $\mathfrak{o}$ is clearly not inverse to $\mathfrak{a}$.  Thus, we must look elsewhere.  In the next section, we do this.\medskip

Nonetheless, before we seek the inverse interpretation, we note in passing that we cannot adapt $\mathfrak{o}$ to give an interpretation of $\ID$ in $\EA^*$.  This is a consequence of the following fact:  in $\EA^*$, we cannot prove that the class of von Neumann ordinals is closed under ordinal addition, let alone multiplication and exponentiation.  This, in turn, is a consequence of the following theorem.
\begin{theorem}\label{gblfnboundinglemma} Suppose $\Phi$ is a bounded quantifier formula of $\EA^*$.  Then 
\[
\EA^* \vdash (\forall \vec{x})(\exists y)\Phi(\vec{x}, y)
\]
if, and only if, there is a natural number $\mathbf{k}$ such that
\[
\EA^* \vdash (\forall \vec{x})(\exists y \in \P^\mathbf{k}(\R(\{x_1, ..., x_n\})))\Phi(\vec{x}, y),
\]
where $\P^\mathbf{k}(x) = \underbrace{\P(\P( ... \P}_{\mathbf{k}}(x) ...))$.
\end{theorem}
\emph{Proof sketch. } This is proved in two steps. First, for bounded quantifier $\Phi$, we show that $\EA^* \vdash (\forall \vec{x})(\exists y)\Phi(\vec{x}, y)$ if, and only if, there is a term $t$ of $\EA^*$ such that $\EA^* \vdash (\forall \vec{x})(\exists y \in t(\vec{x}))\Phi(\vec{x}, y)$.  Clearly, this is analogous to Parikh's celebrated result concerning $\ID$ (Theorem 4.4 of \cite{parikh1971eaf}), and may be proved using a similar compactness argument.  This is possible in part because Separation is restricted to bounded quantifier formulae; a single unbounded quantifier instance would render the theorem false.  Second, we show, by induction on the construction of terms in $\EA^*$, that, for any term $t$ of $\EA^*$, there is $\mathbf{k}$ such that $\EA^* \vdash \forall \vec{x}(t(\vec{x}) \in \P^\mathbf{k}(\R(\{x_1, ..., x_n\})))$. \hfill $\Box$\medskip

\noindent As we will see in Section \ref{IDEA}, this result also entails that Kaye and Wong's inverse to the Ackermann interpretation cannot be defined from $\ID$ to $\EA^*$.  

However, although we cannot define an \emph{ordinal} interpretation of $\ID$ in $\EA^*$, we can define a \emph{cardinal} interpretation: see \cite{homolka1983ea}.  To state this, we need some notation:
\begin{itemize}
\item  $x \leq_\mathrm{c} y$ iff there is an injection from $x$ into $y$
\item $x \simeq_\mathrm{c} y$ iff $x \leq_\mathrm{c} y$ and $y \leq_\mathrm{c} x$
\item $x <_\mathrm{c} y$ iff $x \leq_\mathrm{c} y$ but $y \not \leq_\mathrm{c} x$.
\end{itemize}
Let $\mathfrak{c} : \ID \rightarrow \EA^*$ be the interpretation mapping defined as follows:
\begin{eqnarray*}
\mathrm{Dom}(x)^\mathfrak{c} & \mbox{ is } & \mathrm{Dom}(x) \\
(x = y)^\mathfrak{c} & \mbox{ is } & (x \simeq_\mathrm{c} y) \\
0^\mathfrak{c} & \mbox{ is } & \0 \\
(x < y)^\mathfrak{c}  & \mbox{ is } & (x <_\mathrm{c} y) \\
(S(x))^\mathfrak{c} & \mbox{ is } & x \cup \{x\} \\
(x + y)^\mathfrak{c} & \mbox{ is } & (x \times \{\0\}) \cup (y \times \{\{\0\}\}) \\
(x \cdot y)^\mathfrak{c} & \mbox{ is } & x \times y \\
(\mathrm{Exp}(x, y))^\mathfrak{c} & \mbox{ is } & \{f : y \rightarrow x\}
\end{eqnarray*}
\begin{theorem}
$\mathfrak{c}$ defines an interpretation of $\ID$ in $\EA^*$.
\end{theorem}
Under this interpretation, the bounded induction axioms of $\ID$ follow from $\in$-induction for bounded quantifier formulae in $\EA^*$:  see Theorem 8.3.3 of \cite{mayberry2000fmt}.  Again, however, it is clear that $\mathfrak{c}$ and $\mathfrak{a}$ are not inverses of each other.

\section{The inverse to the Ackermann interpretation}\label{IDEA}

To define the inverse to the Ackermann interpretation of $\ZFinf^*$ in $\mathsf{PA}$, Kaye and Wong exploit a function $\mathfrak{p} : V \rightarrow \mathrm{Ord}$, which takes each set to its `Ackermann code' in the von Neumann ordinals.  That is, $\mathfrak{p}$ satisfies the following $\in$-recursive definition,
\[
\mathfrak{p}(x) = \sum_{y \in x} 2^{\mathfrak{p}(y)}
\]
where the bounded sum and exponentiation operation on the right-hand side are ordinal bounded sum and ordinal exponentiation respectively.  With this in hand, they define $\mathfrak{b} : \mathsf{PA} \rightarrow \mbox{ZF--Inf}^*$ as follows:
\begin{eqnarray*}
\mathrm{Dom}(x)^\mathfrak{b} & \mbox{ is } & \mathrm{Dom}(x) \\
(x = y)^\mathfrak{b} & \mbox{ is } & x = y \\
(x < y)^\mathfrak{b} & \mbox{ is } & \mathfrak{p}(x) < \mathfrak{p}(y) \\
(x + y = z)^\mathfrak{b} & \mbox{ is } & \mathfrak{p}(x) + \mathfrak{p}(y) = \mathfrak{p}(z) \\
(x \cdot y = z)^\mathfrak{b} & \mbox{ is } & \mathfrak{p}(x) \times \mathfrak{p}(y) = \mathfrak{p}(z)
\end{eqnarray*}
where the relations in the final three lines on the right-hand side are relations on the ordinals.

In Section \ref{ordcardinterps}, I noted that the von Neumann ordinals are not closed under addition in $\EA^*$ and I remarked that this precludes the usual ordinal interpretation of $\ID$ in $\EA^*$.  Here again it prevents an interpretation.  Clearly, we cannot define Kaye and Wong's function $\mathfrak{p}$, nor \emph{a fortiori} their interpretation $\mathfrak{b}$.  Thus, we must be more resourceful.\medskip

Essentially, Kaye and Wong's inverse interpretation $\mathfrak{b}$ exploits two facts: (i) the von Neumann ordinals provide a model of $\mathsf{PA}$; and (ii) there is a bijection between the universe and that model that takes a set to its `Ackermann code' in the model.  We cannot adapt their construction because, as we have seen, in $\EA^*$, the von Neumann ordinals do not provide a model of $\ID$.  

However, we can adapt their strategy.  I will define a class of sets in $\EA^*$ with the following two properties:  (i) it provides a model of $\ID$; and (ii) there is a bijection between the universe and that model that takes a set to its `Ackermann code' in the model.

Indeed, the class is $V$, the class of all sets.  And the bijection is simply the identity mapping.  That is, I will define a set $0_\mathrm{a}$, a relation $<_\mathrm{a}$, and functions $S_\mathrm{a}$, $+_\mathrm{a}$, $\times_\mathrm{a}$, and $\mathrm{Exp}_\mathrm{a}$ such that $\<V, 0_\mathrm{a}, <_\mathrm{a}, S_\mathrm{a}, +_\mathrm{a}, \times_\mathrm{a}, \mathrm{Exp}_\mathrm{a}\> \models \ID$.  Then I will show that each set is its own Ackermann code, when considered as an element in this model.  This will give rise to the following natural interpretation $\mathfrak{d} : \ID \rightarrow \EA^*$, which is inverse to $\mathfrak{a}$:
\begin{eqnarray*}
\mathrm{Dom}(x)^\mathfrak{d} & \mbox{ is } & \mathrm{Dom}(x) \\
(x = y)^\mathfrak{d} & \mbox{ is } & x = y \\
0^\mathfrak{d} & \mbox{ is } & 0_\mathrm{a} \\
(x < y)^\mathfrak{d} & \mbox{ is } & x <_\mathrm{a} y \\
(S(x))^\mathfrak{d} & \mbox{ is } & S_\mathrm{a}(x) \\
(x + y)^\mathfrak{d} & \mbox{ is } & x +_\mathrm{a} y \\
(x \cdot y)^\mathfrak{d} & \mbox{ is } & x \times_\mathrm{a} y \\
(\mathrm{Exp}(x, y))^\mathfrak{d} & \mbox{ is } & \mathrm{Exp}_\mathrm{a}(x, y)
\end{eqnarray*}

The definitions of $0_\mathrm{a}$, $<_\mathrm{a}$, $S_\mathrm{a}$, $+_\mathrm{a}$, $\times_\mathrm{a}$, and $\mathrm{Exp}_\mathrm{a}$ depend on a function that takes each level of the cumulative hierarchy $V_n$ to a linear ordering of $V_n$.  To define this function, we need to introduce the notion of a \emph{lexicographical ordering}.  First, notation:  given a linear ordering $L = [x_0, ..., x_n]$, let $\Field(L) =_{df.} \{x_0, ..., x_n\}$.
\begin{definition}
Given a linear ordering $L$, define the lexicographical ordering, $\mathrm{Lex}(L)$, of the power set of $\Field(L)$ as follows:  given $X, Y \subseteq \Field(L)$,
\[
X <_{\mathrm{Lex}(L)} Y \mbox{ iff the $L$-greatest element of $X \bigtriangleup Y$ is in $Y$}
\]
where $X \bigtriangleup Y$ is the symmetric difference of $X$ and $Y$.
\end{definition}

Now suppose $V_n$ is a level of the cumulative hierarchy.  That is, there is a linear ordering $[V_0, ..., V_n]$ such that $V_0, ..., V_n \subseteq V_n$, $V_0 = \0$, and $V_{k+1} = \P(V_k)$ for $k = 0$, ..., $n-1$.  Then we define a local function
\[
\Ack : \{V_0, ..., V_n\} \rightarrow \{L \subseteq V_n \times V_n: L \mbox{ is a linear ordering}\}
\]
We define $\Ack$ by recursion along $[V_0, ..., V_n]$ as follows:
\begin{eqnarray*}
\Ack(V_0) & = & [\, ] \mbox{\ \ (the empty ordering)} \\
\Ack(V_{k+1}) & = & \mathrm{Lex}(\Ack(V_k))
\end{eqnarray*}
In $\EA^*$, we can prove that this recursion is effective---that is, we can prove that there is such a local function $\Ack$.  The reason is that a set containing all the values taken by $\Ack$ can be specified prior to carrying out the recursion:  the set is the set of those linear orderings whose fields are subsets of $V_n$.  Thus, it is an instance of \emph{definition by limited recursion} in Mayberry's terminology: see Theorem 9.2.2 of \cite{mayberry2000fmt}.  In $\EA^*$, recursions may be carried out if it is possible to specify a set containing the range of the recursively defined function \emph{prior to defining the function}.  Recursions in which this is not possible are not necessarily effective in $\EA^*$.

It is easy to show that $\Ack(V_0) \subseteq_* \Ack(V_1) \subseteq_* ... \subseteq_* \Ack(V_{n-1}) \subseteq_* \Ack(V_n)$.  Thus, together with the Weak Hierarchy Principle, this construction induces an order on the universe of sets:
\[
x <_\mathrm{a} y \mbox{ iff } \R(x) \subseteq \R(y) \wedge x <_{\Ack(\R(y))} y
\]
(Recall that $\R(x)$ is the first level of the cumulative hierarchy at which $x$ occurs.)  I claim that the universe of sets, ordered in this way, provides a model of $\ID$.  I describe this model precisely now.\medskip

Let $0_\mathrm{a} = \0$.  Let $S_\mathrm{a}(x)$ be the element of the ordering $\Ack(\P(\R(x)))$ that follows immediately after $x$.  (Since all sets and thus all linear orderings are finite, every linear ordering has endpoints and immediate successors and predecessors.)

To define addition, multiplication, and exponentiation, we require a little notation: Given a linear ordering $L$ and $x, y \in \Field(L)$, let $[x, ..., y]_L$ denote the segment of $L$ between $x$ and $y$ inclusive.

Now, without loss of generality, suppose $x <_\mathrm{a} y$.  So $x, y \in \P(\R(y))$. Then, since $\Ack(\P(\R(y))$ is an ordering of $\P(\R(y))$, which is a level of the cumulative hierarchy, it follows that $x, y \in \Field(\Ack(\P(\R(y))))$. Then let $x+_\mathrm{a} y$ be the unique $z \in \Ack(\P(\R(y)))$ such that\medskip

$\Field\left ([\{\0\}, ..., x]_{\Ack(\P(\R(y)))}\right ) +_\mathrm{c} \Field\left ([\{\0\}, ..., y]_{\Ack(\P(\R(y)))}\right )$

\hspace{60mm} $\simeq_\mathrm{c} \Field\left ([\{\0\}, ..., z]_{\Ack(\P(\R(y)))}\right )$\medskip

\noindent In the proof of Theorem \ref{numerals}, it will become clear why we must begin with $\{\0\}$ rather than with $\0$:  in short, it avoids a `bug by one' problem.

Define $x \times_\mathrm{a} y$ and $\mathrm{Exp}_\mathrm{a}(x, y)$ similarly.\medskip

This completes our definition of the interpretation mapping $\mathfrak{d} : \ID \rightarrow \EA^*$.  The following easy theorem establishes that $\mathfrak{d}$ defines an interpretation of $\ID$ in $\EA^*$.
\begin{theorem}
$\<V, 0_\mathrm{a}, <_\mathrm{a}, S_\mathrm{a}, +_\mathrm{a}, \times_\mathrm{a}, \mathrm{Exp}_\mathrm{a}\> \models \ID$.
\end{theorem}
We now turn to the problem of showing that $\mathfrak{a}$ and $\mathfrak{d}$ are inverses.  It suffices to prove the following theorem:
\begin{theorem}\label{numerals}
\[
\EA^* \vdash x \in y \leftrightarrow (\mbox{the $x^\mathrm{th}$ bit of $y$ is $1$})^\mathfrak{d}
\]
\end{theorem}
\emph{Proof. } We prove this indirectly.  First, we define a function that assigns to each set a binary numeral:  Given a set $x$, let $\Num(x)$ be the sequence (or binary numeral) $\<s_0, ..., s_n\>_{[x_0, ..., x_n = x]_{\Ack(\R(x))}}$ where
\[
s_i = \left \{ \begin{array}{ll}
1 & \mbox{if } x_i \in x \\
0 & \mbox{if } x_i \not \in x
\end{array}
\right.
\]
Then we note that it follows easily from the definition of lexicographical orderings and $S_\mathrm{a}(x)$ that, if $\Num(x)$ is \medskip

$\<1, 1, ..., 1, 0, s_k, ..., s_n\>_{[x_0, ..., x_k, ..., x_n = x]_{\Ack(\R(x))}}$\medskip

\noindent then\medskip

$\<0, 0, ..., 0, 1, s_k, ..., s_n, 0\>_{[x_0, ..., x_k, ..., x_n = x, x_{n+1} = S_\mathrm{a}(x)]_{\Ack(\R(S_\mathrm{a}(x)))}}$\medskip

\noindent is  $\Num(S_\mathrm{a}(x))$.  Thus, given $x$ and $[x_1, ..., x_n = x]_{\Ack(\R(x))}$, the linear ordering
\[
[\Num(x_1), ..., \Num(x_n)]
\]
contains all binary numerals between $\Num(x_1) = \<1, 0\>_{[x_0, x_1]}$ and $\Num(x_n)$ inclusive.  And, if
\[
\Num(x) = \<s_0, ..., s_n\>_{[x_0, ..., x_n]}
\]
then there are
\[
s_02^0 + \cdots + s_n2^n
\]
such numerals.  Thus,
\begin{eqnarray*}
& & \Field\left ([x_1, ..., x_n = x]_{\Ack(\R(x))}\right ) \\
& \simeq_\mathrm{c} & \Field\left ([\Num(x_1), ..., \Num(x_n)]_{\Ack(\R(x))} \right) \\
& \simeq_\mathrm{c} & s_02^0 + \cdots + s_n2^n
\end{eqnarray*}
Now, suppose $x <_\mathrm{a} y$ and $[x_0, ..., x_m = x, ..., x_n = y]_{\Ack(\R(y))} \subseteq_* \Ack(\R(y))$.  And suppose
\[
\Num(y) = \<s_0, ..., s_m, ..., s_n\>_{[x_0, ..., x_m = x, ..., x_n = y]_{\Ack(\R(y))}}
\]
Then
\begin{eqnarray*}
& & x \in y \\
& \mbox{ iff } & s_m = 1\\
& \mbox{ iff }&  \mbox{the $\Field([x_1, ..., x_m])^\mathrm{th}$ bit of $\Field([x_1, ..., x_m, ..., x_n])$ is $1$} \\
& \mbox{ iff } & (\mbox{the $x^\mathrm{th}$ bit of $y$ is 1})_\mathrm{a} \\
& \mbox{ iff } & (\mbox{the $x^\mathrm{th}$ bit of $y$ is 1})^\mathfrak{d}
\end{eqnarray*}
as required. \hfill $\Box$
\begin{corollary}
$\EA^*$ and $\ID$ are bi-interpretable. The interpretations $\mathfrak{a}$ and $\mathfrak{d}$ are inverse to each other.
\end{corollary}

\section{Concluding remarks}

As in $\ZFinf^*$, the Axiom of Choice is provable in $\EA^*$:  the proof is an easy application of One Point Extension Induction for bounded quantifier formulae. Thus, its Ackermann translation holds in $\ID$.

A little more interesting is the fact that $\EA^*$ does not prove that for every set there is a finite von Neumann ordinal of the same size.  If $\EA^*$ were to prove this, then it would prove that the von Neumann ordinals are closed under exponentiation, which it does not, by Lemma \ref{gblfnboundinglemma}.  Thus, $\ID$ does not prove the translation of this sentence.

More interestingly still, it is not known whether or not $\ID$ proves the Ackermann translation of the bounded replacement scheme:  that is, for each bounded quantifier formula $\Phi$,
\[
\forall x \exists! y \Phi(x, y) \rightarrow \forall x \exists y \forall z (z \in y \equiv (\exists u \in x)\Phi(u, z))
\]
The translation of each such sentence is provable in $\ID$ with the bounded collection scheme for $\Sigma_1$-formulae, but this is the strongest result known.  The equivalence of $\EA^*$ and $\ID$ opens up a new way to investigate this question.  The following result is the best known in $\EA^*$:
\begin{theorem}
Suppose $\Phi$ is a bounded quantifier formula of $\EA^*$.  Then, if
\[
\EA^* \vdash \forall x \exists! y \Phi(x, y)
\]
then
\[
\EA^* \vdash \forall x \exists y \forall z (z \in y \equiv (\exists u \in x)\Phi(u, z))
\]
\end{theorem}
The proof relies on the Parikh-style result used in the proof of Theorem \ref{gblfnboundinglemma}.

\footnotesize{
\bibliography{interparith}{}
\bibliographystyle{plain}
}

\end{document}